\documentclass[a4,12pt,reqno,fleqn]{amsart}
\usepackage{graphicx}
\usepackage{amscd}
\usepackage{xfrac}
\usepackage{mathtools}
\emergencystretch=\maxdimen
\allowdisplaybreaks[1]
\usepackage{array}
\usepackage{xcolor}
\usepackage[all]{xy}
\makeatletter
\@namedef{subjclassname@2020}{\textup{2020} Mathematics Subject Classification}
\makeatother


\usepackage{amsmath}
\usepackage{amsthm}
\usepackage{amsfonts,amssymb,hyperref,mathrsfs,cleveref,setspace,enumitem,verbatim}
\usepackage[margin=3.5cm]{geometry}

\DeclareMathAlphabet{\mathpzc}{OT1}{pzc}{m}{it}

\newtheorem{thm}{Theorem}[section]

\newtheorem{lem}[thm]{Lemma}
\newtheorem{propn}[thm]{Proposition}
\theoremstyle{definition}
\newtheorem{defn}[thm]{Definition}

\newcommand{\bj}{\perp_{BJ}}

\author[Mohit and R.~Jain]{Mohit and Ranjana Jain}

\address{Mohit, Department of  Mathematics, University of Delhi, Delhi, India}
\email{mohitdhandamaths@gmail.com}

\address{Ranjana Jain, Department of Mathematics, University of Delhi, Delhi, India}
\email{rjain@maths.du.ac.in}
\thanks{Research of the first named author is supported by Savitribai Jyotirao Phule Single Girl Child Fellowship vide F.No. 82-7/2022(SA-III)}
\begin{document}
	\title[Approximate B-J orthogonality preserver on Lebesgue-Bochner spaces]{Approximate Birkhoff-James orthogonality preserver on Lebesgue-Bochner spaces}
	\maketitle
	\textbf{Abstract:} In this article, we examine an approximate version of Koldobsky-Blanco-Turn\v{s}ek theorem (namely, Property P) in the space of vector-valued integrable functions. More precisely, we prove that the Lebesgue-Bochner spaces $L^p(\mu,X),\;(1\leq p<\infty)$, do not have Property P under certain conditions on $\mu$ and the Banach space $X$. 
	\section{introduction}
Birkhoff-James (in short, B-J) orthogonality    plays a vital role in the study of geometry of normed spaces. Koldobsky \cite{koldo}, in 1993, and Blanco and Turn\v{s}ek \cite{blanco}, in 2006, utilized the concept of B-J orthogonality to characterize the isometries on Banach spaces (for real and complex, respectively). In particular, they proved that  any linear operator on a Banach space that preserves the B-J orthogonality must be a scalar multiple of an isometry. This study was further carried over in the context of geometric non-linear characterizations of certain Banach spaces \cite{tana1,tana2,tana3}. Apart from this, several refinements of Koldobsky-Blanco-Turn\v{s}ek theorem have also been studied by researchers (see,  \cite{manna, sain}). 

	
	 An element $x$ in a normed space $X$ over $\mathbb{F} \ (\mathbb{R} \;\text{or}\; \mathbb{C})$ is said to be {\it Birkhoff-James orthogonal} to $y\in X$ if $\|x+\alpha y\|\geq\|x\|$ for all $\alpha\in\mathbb{F}$. Later, Chmieli\'{n}ski\cite{jacek} introduced an approximate version of B-J orthogonality. For $\epsilon\in(0,1)$, an element $x\in X$ is said to be {\it $\epsilon$-approximate B-J orthogonal} to $y\in X$ (written as, $x\perp^{\epsilon}_{BJ}y$) if the following inequality holds: $$\|x+\alpha y\|^2\geq\|x\|^2-2\epsilon\|x\|\|\alpha y\|,\;\forall\;\alpha\in\mathbb{F}.$$ 
	 
	
	 Recently, Mandal et.al.\cite{mandal} introduced an approximate version of preservation of orthogonality  by a linear operator defined on a Banach space. More explicitly, a Banach space $X$ is said to have {\it Property P} if any bounded linear operator defined on $X$ which preserves approximate B-J orthogonality is a scalar multiple of an isometry. Mandal et.al. exhibited that the space $\ell^{n}_{\infty}$ for $n>2$ and a two-dimensional Banach space whose unit sphere is a regular polygon possess Property P. They also established that $\ell^p,\;1\leq p<\infty$, the space of absolutely $p$-summable real valued sequences does not have Property P. 
	
Our aim in this article is to explore Property P for the Lebesgue-Bochner spaces $L^p(\mu,X),\;1\leq p<\infty$, the spaces of  vector-valued integrable functions. We provide  descriptions of a support functionals at  non-zero elements of $\ell^1(X)$, the space of vector-valued absolutely summable sequences, and of $L^1(\mu,X)$. Further, we define a semi-inner product on $L^p(\mu,X),\;1<p<\infty$, where $\mu$ is a complete positive measure and  norm of $X$ is Fr\'{e}chet differentiable. These results help us to establish that  the spaces $\ell^1(X)$ and $L^p(\mu,X)$, $1\leq p<\infty$, do not have Property $P$.

	\section{Preliminaries}
 In this section, we give few definitions and results which will be needed for further investigation. 
 
	\begin{defn}\cite[Definition 1.2]{mandal}
		A Banach space $X$ over $\mathbb{F}$ is said to have {\it Property P} if there exists an $\epsilon\in (0,1)$ such that any bounded linear operator  $T:X \to X$ which preserves $\epsilon$-approximate orthogonality (that is, for any $x,\;y\in X$,\; $x\bj y$ $\implies$ $T(x)\perp^{\epsilon}_{BJ}T(y)$ ) must be a scalar multiple of an isometry.
	\end{defn}
Clearly, Property $P$ can be seen as an approximate version of Koldobsky-Blanco-Turn\v{s}ek theorem.

 It is evident that James criteria is an important tool for B-J orthogonality. Next, we state an analogous result which characterizes $\epsilon$-approximate B-J orthogonality in Banach spaces.
\begin{thm}\cite[Theorem 2.3]{stypula}\label{approxjames}
	Let $X$ be a real Banach space and $x,y\in X$. Then, for $\epsilon\in[0,1)$,\; $x\perp^{\epsilon}_{BJ}y$ if and only if there exists a bounded linear functional satisfying $\|T\|=1,\;T(x)=\|x\|$ and $|T(y)|\leq\epsilon\|y\|$. Moreover, reverse implication is also true for complex Banach spaces.
\end{thm}

With the aim of generalizing inner product in normed spaces, Giles\cite{giles} introduced the concept of semi-inner product in normed spaces.
\begin{defn}\cite{giles}
	 A mapping $[\cdot]:X\times X\rightarrow \mathbb{F}$ is said to be a {\it semi-inner product} on $X$ if for any arbitrary $x, y,z\in X$ and $\alpha, \beta\in \mathbb{F}$, the following hold true:
	\begin{enumerate}
		\item $[\alpha x+\beta y,z]=\alpha[x,z]+\beta[y,z]$,
		\item $[x,\alpha y]=\bar{\alpha}[x,y]$,
		\item $|[x,y]|\leq\|x\|\|y\|$,
		\item $[x,x]=\|x\|^2$.
	\end{enumerate}
\end{defn}
\begin{defn}
		For a measure space $(S,\mu)$ and a Banach space $X$ (real or complex), the {\it Lebesgue-Bochner space} $L^{p}(\mu, X),\;1\leq p<\infty$, is defined as
	$$L^{p}(\mu, X):=\{f:S\rightarrow X|\ f \ \text{is strongly measurable and} \ \;\int\limits\limits\limits_{S}\|f(s)\|^p\,ds <\infty\},
	$$
	where almost everywhere equal functions are identified.
\end{defn}
\begin{defn}
		The norm function $\| \cdot \|$ on $X$ is said to be {\it Fr\'{e}chet differentiable} at $0\neq x\in X$ if there exists $f\in X^{*}$ such that $$\lim_{h \to 0}\frac{\big|\|x+h\|-\|x\|-f(h)\big|}{\|h\|}=0.$$ 
\end{defn}
The norm function is called {\it Fr\'{e}chet differentiable} if it is Fr\'{e}chet differentiable at every non-zero point.

		\begin{thm}\cite{ion}\label{duall1}
		The space $\ell^1(X)^*$ is isometrically isomorphic to $\ell^{\infty}(X^*)$ via the mapping $\theta:\ell^{\infty}(X^*)\rightarrow\ell^1(X)^*$ defined as $(f_n)\mapsto \theta((f_n))$, where $\theta((f_n))((y_n))=\sum\limits_{n\in \mathbb{N}}f_n(y_n)$, for all $(y_n)\in \ell^1(X)$. 
	\end{thm}
\begin{thm}\cite[Theorem 3]{chin}\label{dualL1}
	Let $(S,\Sigma,\mu)$ be a $\sigma$-finite, positive measure space and let $X$ be a Banach space (real or complex) such that $X^*$ has the wide Radon-Nikodym property with respect to $\mu$. Then $L^1(\mu,X)^*$ is isometrically isomorphic to $L^{\infty}(\mu,X^*)$ via the mapping $L^{\infty}(\mu,X^*) \ni G\rightarrow \theta(G) \in L^1(\mu,X)^*$, where $\theta(G)(f)=\int\limits_{S}G(s)(f(s))\,ds$, for all $f\in L^1(\mu,X)$.
\end{thm}

Throughout the article, $\mu$ denotes a complete positive measure and $X$ denotes a real Banach space unless specified. For $f\in L^p(\mu,X),\;1\leq p<\infty$, we denote the set $\{s\in S:\;f(s)=0\}$ by $Z(f)$ and for a non-zero element $x$ in a Banach space $X$, $J(x)=\{f\in X^*: \|f\|=1,\;f(x)=\|x\|\}$ denotes the collection of support functionals at $x$.

	\section{main results}
	We start the section by characterizing the support functional corresponding to a non-zero element in $\ell^{1}(X)$.
	\begin{lem}\label{supportl1}
		 Let $x=(x_n)\in \ell^1(X)$ be a non-zero element. Then, $T\in J(x)$ if and only if $T((y_n))=\sum\limits_{n\in \mathbb{N}}F_n(y_n)$ for every $(y_n)\in \ell^1(X)$, where $(F_n)\in \ell^{\infty}(X^*)$ with $\|F_n\|\leq1$, for all $n\in\mathbb{N}$ and $F_n\in J(x_n)$, for all $n\in \mathbb{N}$ such that $x_n\neq0$.
		\end{lem}
	\begin{proof}
 Suppose that $T\in J(x)$. By \Cref{duall1}, there exists $(H_n)\in \ell^{\infty}(X^*)$ such that $\|(H_n)\|_{\infty}=1$ and $T((y_n))=\sum\limits_{n\in \mathbb{N}}H_n(y_n)$. Now, $T(x)=\|x\|$ which gives that $\sum\limits_{\{n\in\mathbb{N}: x_n\neq0\}}(\|x_n\|-H_n(x_n))=0$. Since $\|H_n\|\leq1$, for all $n\in\mathbb{N}$ and $X$ is a real Banach space, therefore, $H_n(x_n)\leq\|x_n\|$, for all $n\in\mathbb{N}$. Thus, $\sum\limits_{\{n\in\mathbb{N}:\ x_n\neq0\}}(\|x_n\|-H_n(x_n))=0$ implies that $H_n(x_n)=\|x_n\|$ whenever $x_n\neq0$. Therefore, for $n\in \mathbb{N}$ such that $x_n\neq0$, we have,
	$$	\|H_n\|=\sup\limits_{x\in X\setminus{\{0\}}}\frac{|H_n(x)|}{\|x\|}\geq\frac{H_n(x_n)}{\|x_n\|}=1.$$
	This gives that $\|H_n\|=1$ and hence $H_n\in J(x_n)$, for all $n\in\mathbb{N}$ such that $x_n\neq0$. This shows that $T$ corresponds to a sequence $(H_n)\in \ell^{\infty}(X^*)$ such that $H_n\in J(x_n)$, if $x_n\neq0$ and $\|H_n\|\leq1$, if $x_n=0$. 
	
	Conversely,	since $x=(x_n)\in \ell^1(X)$ is a non-zero element and $\|F_{n}\|=1$, whenever $x_n\neq0$, therefore, $\|(F_n)\|_{\infty}=1$.  Thus, by \Cref{duall1}, $T$ is a norm-one linear functional and $T((x_n))=\|(x_n)\|$. This completes the proof.
	\end{proof}
By utilizing the above result, we prove that the approximate version of  Koldobsky-Blanco-Turn$\check{s}$ek theorem does not hold in $\ell^{1}(X)$.
\begin{thm}
	The sequence space $\ell^1(X)$ does not have Property P.
\end{thm}
\begin{proof}
	 We prove that for each $\epsilon\in(0,1)$, there exists $U_{\epsilon}\in B(\ell^1(X))$, the space of bounded operators defined on $\ell^1(X)$, which preserve the $\epsilon$-approximate B-J orthogonality but is not a scalar multiple of an isometry. For $\epsilon\in(0,1)$, define $U_{\epsilon}:\ell^1(X)\rightarrow \ell^1(X)$ as $U_{\epsilon}((x_n))=((1-\epsilon)x_1,\;x_2,...)$. Clearly, $U_{\epsilon}\in B(\ell^1(X))$ and is not a scalar multiple of an isometry. 
	
	Now, we claim that $U_{\epsilon}$ preserves the $\epsilon$-approximate B-J orthogonality. Let $x=(x_n),\;y=(y_n)\in \ell^1(X)$ be any arbitrary such that $x\bj y$. By James criteria, there exists $T\in J(x)$ such that $T(y)=0$. By \Cref{supportl1}, $T:\ell^1(X)\rightarrow \mathbb{R}$ can be represented as $T((z_n))=\sum\limits_{n\in\mathbb{N}}F_n(z_n)$, where $(F_n)\in \ell^{\infty}(X^*)$ with $F_n\in J(x_n)$, for all $n$ such that $x_n\neq0$.
 Clearly, $T(U_{\epsilon}(x))=\|U_{\epsilon}(x)\|$. To prove the claim, by \cite[Theorem 2.3]{stypula}, it is sufficient to show $|T(U_{\epsilon}(y))|\leq\epsilon\|U_{\epsilon}(y)\|$. For this, consider, 
 $$|T(U_{\epsilon}(y))|=\bigg|(1-\epsilon)F_1(y_1)+\sum\limits_{n\geq2}F_n(y_n)\bigg|.$$ Now, $T(y)=0$ gives that $F_1(y_1)=-\sum\limits_{n\geq2}F_n(y_n)$. Thus, 
 \begin{align*}
 	|T(U_{\epsilon}(y))|&=\bigg|-(1-\epsilon)\sum\limits_{n\geq2}F_n(y_n)+\sum\limits_{n\geq2}F_n(y_n)\bigg|\\
 	& \leq\epsilon\sum\limits_{n\geq2}\|y_n\|\\
 	&\leq\epsilon\bigg|(1-\epsilon)\|y_1\|+\sum\limits_{n\geq2}\|y_n\|\bigg|\\ &=\epsilon\|U_{\epsilon}(y)\|.
 \end{align*}
Thus, $U_{\epsilon}(x)\perp^{\epsilon}_{BJ}U_{\epsilon}(y)$, and the claim is proved.
\end{proof}
One  may observe that in order to check Property P, the information about the support functional is very crucial. Thus, in order to move to the space $L^1(\mu,X)$, we first derive a useful form of the support functional at a given non-zero element of $L^1(\mu,X)$. 
\begin{lem}\label{supportL1}
		Let $(S,\mu)$ be a $\sigma$-finite measure space and $X$ be a Banach space such that $X^*$ has the wide Radon Nikodym property with respect to $\mu$. Let $f\in L^1(\mu,X)$ be a non-zero element. If $T\in L^1(\mu,X)^*$ is a support functional for $f$, then $T:L^1(\mu,X)\rightarrow \mathbb{R}$ is of the form $T(g)=\int\limits_{S}G(s)(g(s))\,ds$, for all $g\in L^1(\mu,X)$, where $G\in L^{\infty}(\mu,X^*)$ with $G(s)\in J(f(s))$ for a.e. $s\in Z(f)^c$.
\end{lem}
\begin{proof}
	Since $T$ is a support functional at $f$, by \Cref{dualL1}, there exists $G\in L^{\infty}(\mu,X^*)$ such that $T(g)=\int\limits_{S}G(s)(g(s))\,ds$, for all $g\in L^1(\mu,X)$.  Since $T(f)=\|f\|$, we have
	 \begin{equation}\label{eq1sup}
	\int\limits_{Z(f)^c}(\|f(s)\|-G(s)(f(s)))\,ds=0.
	\end{equation}
	 Also, $\|T\|=1$ gives that $\|G\|_{\infty}=\underset{s\in S}{\text{esssup}}\|G(s)\|=1$, that is, $\|G(s)\|\leq1$ for a.e. $s\in S$. Since $\mu$ is complete, we get $|G(s)(f(s))|\leq\|f(s)\|$ for a.e. $s\in Z(f)^c$. Thus, the function $s\mapsto \|f(s)\|-G(s)(f(s))$ is a non-negative function for a.e. $s\in Z(f)^c$,\; $X$ being real Banach space. Hence, by \Cref{eq1sup}, $G(s)(f(s))=\|f(s)\|$ for a.e. $s\in Z(f)^c$.  Hence, $G(s)\in J(f(s))$ for a.e. $s\in Z(f)^c$.
\end{proof}
\begin{thm}
		Let $(S,\Sigma,\mu)$ be a $\sigma$-finite measure space with atleast two disjoint measurable sets of finite positive measures. Let $X$ be a Banach space such that $X^*$ has the wide Radon Nikodym property with respect to $\mu$. Then the space $L^1(\mu,X)$ does not have Property P.
\end{thm}
	\begin{proof}
		Let $A,\;B\in \Sigma$ be such that $0<\mu(A),\;\mu(B)<\infty$, and $A\cap B\neq\emptyset$. For an arbitrary $\epsilon\in(0,1)$, define $U_\epsilon:L^1(\mu,X)\rightarrow L^1(\mu,X)$ as $U_\epsilon(f)=(1-\epsilon)f\chi_{A}+f\chi_{S\setminus A}$, where $\chi_A$ denotes the characteristic function corresponding to the set $A$. Clearly, $U_\epsilon\in B(L^1(\mu,X))$. 
		
		Observe that $U_\epsilon$ is not a scalar multiple of an isometry. For this, let $\alpha\in \mathbb{R}$ be any arbitrary and let $x_0\in X$ be a norm one element. Define a function $h_\alpha: S\rightarrow X$ as $h_{\alpha}(s)=\chi_{A}(s)x_0+\alpha\chi_{B}(s)x_0$, then $h_\alpha\in L^1(\mu,X)$ with $\|h_\alpha\|=\mu(A)+|\alpha|\mu(B)$. Also, $U_\epsilon(h_\alpha)=(1-\epsilon)\chi_{A}x_0+\alpha\chi_{B}x_0$, and $\|U_\epsilon(h_\alpha)\|=(1-\epsilon)\mu(A)+|\alpha|\mu(B)$. Let, if possible, $U_\epsilon$ is a scalar multiple of an isometry, then there exists a $c>0$ such that $\|U_\epsilon(f)\|=c\|f\|$, for all $f\in L^1(\mu,X)$. In particular, $\|U_\epsilon(h_\alpha)\|=c\|h_\alpha\|$, and this gives that $c=\frac{(1-\epsilon)\mu(A)+|\alpha|\mu(B)}{\mu(A)+|\alpha|\mu(B)}$, since $c$ depends on a scalar $\alpha$ which leads to a contradiction. Thus, $U_\epsilon$ is not a scalar multiple of an isometry.
		
		Now, we claim that $U_\epsilon$ preserves the $\epsilon$-approximate B-J orthogonality. Let $f,\;g\in L^1(\mu,X)$ be such that $f\bj g$. Then by James criteria, there exists $T\in L^1(\mu,X)^*$ satisfying $\|T\|=1,\; T(f)=\|f\|$ and $T(g)=0$. Now, by \Cref{dualL1}, there exists $G\in L^\infty(\mu,X^*)$ such that $T(h)=\int\limits_{S}G(s)(h(s))\,ds$, for all $h\in L^1(\mu,X)$.

		Observe that the set $K=\{s\in Z(f)^c: G(s)\in J(f(s))\}$ is measurable. To see this, let us write $K=K_{1}\cap K_2$, where $K_1=\{s\in Z(f)^c: G(s)(f(s))=\|f(s)\|\}$ and $K_2=\{s\in Z(f)^c:\|G(s)\|\leq1\}$. The map   $s\mapsto G(s)(f(s))$ is a measurable on $S$ and thus, on $Z(f)^c$. Also,  $s\mapsto\frac{1}{\|f(s)\|}$ is measurable on $Z(f)^c$. Together, we get that the map $\psi:\;s\mapsto\frac{1}{\|f(s)\|}G(s)(f(s))$ is measurable on $Z(f)^c$. Thus, the set $K_1=\psi^{-1}(\{1\})$ is a measurable set. Further, $G:S\rightarrow X^*$ is strongly measurable and $\mu$ is a complete measure on $S$, thus, $G$ is  measurable \cite[Lemma 10.2]{lche}. So, $s\mapsto \|G(s)\|$ is a measurable on $Z(f)^c$ which gives the measurability of $K_2$.

		In view of \Cref{supportL1}, we have $G(s)\in J(f(s))$ for a.e. $s\in Z(f)^c$. Therefore, $\mu(Z(f)^c \setminus K) = 0$, and we have
		\begin{align*}
			T(U_\epsilon(f))&=\int\limits_{A\cap Z(f)^c}(1-\epsilon)G(s)(f(s))\,ds+\int\limits_{(S\setminus A)\cap Z(f)^c}G(s)(f(s))\,ds\\
			&=\int\limits_{A\cap K}(1-\epsilon)G(s)(f(s))\,ds+\int\limits_{(S\setminus A)\cap K}G(s)(f(s))\,ds\\
			&=\int\limits_{A\cap K}(1-\epsilon)\|f(s)\|\,ds+\int\limits_{(S\setminus A)\cap K}\|f(s)\|\,ds\\
				&=\int\limits_{A\cap Z(f)^c}(1-\epsilon)\|f(s)\|\,ds+\int\limits_{(S\setminus A)\cap Z(f)^c}\|f(s)\|\,ds\\
				&=\int\limits_{A}(1-\epsilon)\|f(s)\|\,ds+\int\limits_{S\setminus A}\|f(s)\|\,ds\\
				&=\int\limits_{S}\|U_\epsilon(f)(s)\|\,ds=\|U_\epsilon(f)\|.
		\end{align*}
		Also, $T(g)=\int\limits_{A}G(s)(g(s))\,ds+\int\limits_{S\setminus A}G(s)(g(s))\,ds=0$, therefore, 
		\begin{align*}
			T(U_\epsilon(g))&=(1-\epsilon)\int\limits_{A}G(s)(g(s))\,ds+\int\limits_{S\setminus A}G(s)(g(s))\,ds\\
			&=-(1-\epsilon)\int\limits_{S\setminus A}G(s)(g(s))\,ds+\int\limits_{S\setminus A}G(s)(g(s))\,ds\\
			&=\epsilon\int\limits_{S\setminus A}G(s)(g(s))\,ds.
		\end{align*}
	Now, $\|G(s)\|\leq1$ for a.e. $s\in S$ and $\mu$ is complete, so $\|G(s)\|\leq1$ for a.e. $s\in S\setminus A$. Write $S\setminus A=M\cup M'$, where $M=\{s\in S\setminus A: \|G(s)\|>1\}$ is a zero measure set. Then,
	\begin{align*}
		|T(U_\epsilon(g))|&=\epsilon\Bigg|\int\limits_{S\setminus A}G(s)(g(s))\,ds\Bigg|\\
		&= \epsilon\Bigg|\int\limits_{M'}G(s)(g(s))\,ds\Bigg|\\
		&\leq\epsilon\int\limits_{M'}\|G(s)\|\|(g(s))\|\,ds\\
			&\leq\epsilon\int\limits_{M'}\|g(s)\|\,ds\\
				&\leq\epsilon\Bigg[(1-\epsilon)\int\limits_{A}\|g(s)\|\,ds+\int\limits_{M'}\|g(s)\|\,ds\Bigg]\\
				&=\epsilon\Bigg[(1-\epsilon)\int\limits_{A}\|g(s)\|\,ds+\int\limits_{S\setminus A}\|g(s)\|\,ds\Bigg]\\
				&=\epsilon\|U_\epsilon(g)\|.
	\end{align*}
Hence, again by \Cref{approxjames}, $U_{\epsilon}(f)\perp^{\epsilon}_{BJ}U_{\epsilon}(g)$, and this completes the proof of the claim.
	\end{proof}
Lastly, we prove that the space  $L^p(\mu,X),\;1<p<\infty$ does not satisfy Property P. Before that, we define a semi-inner product on  $L^p(\mu,X),\;1<p<\infty$. Note that if $X$ is a smooth Banach space, then $J(x)$ consists of a single element which we denote by $F_x$. Observe that $F_{\alpha x}=\frac{\overline{\alpha}}{|\alpha|}F_x$, for all $\alpha\in\mathbb{C}$.
\begin{propn}\label{semiLp}
	Let $(S,\;\Sigma,\;\mu)$ be a measure space and let $X$ be a Banach space over $\mathbb{F}$ whose norm is Fr$\acute{e}$chet differentiable. Then, for $1<p<\infty$, the map $[\cdot]:L^p(\mu,X)\rightarrow \mathbb{F}$ defined by 
	\begin{align*}
		[f,g]=
		\begin{cases}
			\frac{1}{\|g\|^{p-2}}\int\limits_{Z(g)^c}\|g(s)\|^{p-1}F_{g(s)}(f(s))\,ds;\;g\neq0,\\
			0\hspace{5.8cm};\;g=0
		\end{cases}
	\end{align*}
is a semi-inner product on $L^p(\mu,X)$.
	\end{propn}
\begin{proof}
	First, we show that the map $[\cdot]$ is well defined. For this, let $0\neq f,\;g\in L^p(\mu,X)$ be arbitrary. Observe that the integrand $s\mapsto\|g(s)\|^{p-1}F_{g(s)}(f(s))$ is a measurable function on $Z(g)^c$. To see this, define $\phi:S\rightarrow \mathbb{F}$ as 
	\begin{align*}
		\phi(s)=
		\begin{cases}
			\|g(s)\|^{p-1}F_{g(s)}(f(s));\;s\in Z(g)^c,\\
			0\hspace{3.4cm};\;\text{otherwise}.
		\end{cases}
	\end{align*} 
	 Let $\{f_n\}$ and $\{g_n\}$ be sequences of simple measurable functions which converge to $f$ and $g$, respectively for a.e. $s\in S$. For each $n\in \mathbb{N}$, let $K_n=\{s\in S: g_n(s)\neq0\}$ and define   $\phi_n:S\rightarrow \mathbb{F}$ as 
	\begin{align*}
		\phi_n(s)=
		\begin{cases}
			\|g_n(s)\|^{p-1}F_{g_n(s)}(f_n(s));\;s\in Z(g)^c\cap K_n,\\
			0\hspace{3.9cm};\;\text{otherwise}.
		\end{cases}
	\end{align*} 
	We claim that $\phi_n(s)$ converges to $\phi(s)$ for a.e. $s\in S$. Let $s\in S$ be such that $\underset{n}\lim\;f_n(s)=f(s)$ and $\underset{n}\lim\; g_n(s)=g(s)$.  If $s\in Z(g)$, then $\phi_{n}(s)=0=\phi(s)$. Thus, $\{\phi_n\}$ is a sequence of simple measurable functions which converges to $\phi$ for a.e. $s\in S$.  If $s\in Z(g)^c$, then $\phi(s)=\|g(s)\|^{p-1}F_{g(s)}(f(s))$ and there exists $n_0\in \mathbb{N}$ such that $g_n(s)\in K_n$ for all $n\geq n_0$ as $\underset{n}\lim\;g_n(s)=g(s)$. Therefore, by \cite[Proposition 3.2]{jain}, $F_{g_n(s)}$ converges to $F_{g(s)}$ and hence, $\underset{n}{\lim}\;\phi_n(s)=\underset{n}{\lim}\;\|g_n(s)\|^{p-1}F_{g_n(s)}(f_n(s))=\phi(s)$. This gives that $\phi$ is measurable, $\mu$ being a complete measure. Hence, $\phi\vert_{Z(g)^c}:s\mapsto\|g(s)\|^{p-1}F_{g(s)}(f(s))$ is measurable.
	
	Now, consider
	\begin{align*}
		|[f,g]|&=	\frac{1}{\|g\|^{p-2}}\Bigg|\int\limits_{Z(g)^c\cap Z(f)^c}\|g(s)\|^{p-1}F_{g(s)}(f(s))\,ds\Bigg|\\
		&\leq 	\frac{1}{\|g\|^{p-2}}\int\limits_{Z(g)^c\cap Z(f)^c}\|g(s)\|^{p-1}\|f(s)\|\,ds.
	\end{align*}
	Since $f,g\in L^p(\mu,X)$, therefore, $s\mapsto \|f(s)\|$ is in $L^p(\mu)$ and $s\mapsto \|g(s)\|^{(p-1)}$ is in $L^q(\mu)$, where $q=\frac{p}{p-1}$. Hence, by H\"{o}lder's inequality, we have
	\begin{align*}
		|[f,g]|&\leq\frac{1}{\|g\|^{p-2}}\Bigg(\int\limits_{Z(g)^c\cap Z(f)^c}\|g(s)\|^p\,ds\Bigg)^{\big(\frac{p-1}{p}\big)}\Bigg(\int\limits_{Z(f)^c\cap Z(g)^c}\|f(s)\|^p\,ds\Bigg)^{\frac{1}{p}}\\
		&\leq \frac{1}{\|g\|^{p-2}}(\|g\|^{p-1})\|f\|\\
		&=\|g\|\|f\|.
	\end{align*}
	
	Clearly, $[\cdot]$ is linear in first component. For any $0\neq\alpha\in\mathbb{F}$ and $0\neq f,\;g\in L^p(\mu,X)$, consider
	\begin{align*}
		[f,\alpha g]&=\frac{1}{\|\alpha g\|^{p-2}}\int\limits_{Z(\alpha g)^c}\|\alpha g(s)\|^{p-1}F_{\alpha g(s)}(f(s))\,ds\\
		&=\frac{1}{|\alpha|^{p-2}\|g\|^{p-2}}\int\limits_{Z(g)^c}|\alpha|^{p-1}\|g(s)\|^{p-1}\frac{\bar{\alpha}}{|\alpha|}F_{g(s)}(f(s))\,ds\\
		&=\bar{\alpha}	\frac{1}{\|g\|^{p-2}}\int\limits_{Z(g)^c}\|g(s)\|^{p-1}F_{g(s)}(f(s))\,ds\\
		&=\bar{\alpha}[f,g].
	\end{align*} 
Thus, $[f,\alpha g]=\bar{\alpha}[f,g]$ and $[f,f]=\|f\|^2$, for all $f,\;g\in L^p(\mu,X)$ and $\alpha\in \mathbb{F}$. Hence, $[\cdot]$ is a semi-inner product on $L^p(\mu,X)$.
\end{proof}
\begin{thm}
Let $(S,\;\Sigma,\;\mu)$ be a measure space with the condition that it has two disjoint measurable sets of finite positive measure and let $X$ be a Banach space over $\mathbb{F}$ whose norm is Fr\'{e}chet differentiable. Then, for $1<p<\infty,\; L^p(\mu,X)$ does not have Property P.	
\end{thm}
\begin{proof}
	  Let $A,\;B$ be two disjoint measurable sets such that $0<\mu(A),\;\mu(B)<\infty$. For $\epsilon\in (0,1)$, define $U_\epsilon:L^p(\mu,X)\rightarrow L^p(\mu,X)$ as $U_\epsilon(f)=f\chi_{A}+\big(1-\frac{\epsilon}{p}\big)f\chi_{S\setminus A}$. Clearly, $U_\epsilon\in B(L^p(\mu,X))$. 
	
	First, we claim that $U_{\epsilon}$ is not a scalar of multiple isometry. For this, let  $0\neq\alpha\in\mathbb{F}$ and $x_0\in X$ be a fixed unit vector. Let if possible, $U_\epsilon$ is a scalar multiple of an isometry. Then there exists a $c>0$ such that $\|T(f)\|=c\|f\|$, for all $f\in L^p(\mu,X)$. In particular, for $h_\alpha=x_0(\chi_{A}+\alpha\chi_{B})\in L^p(\mu,X)$, we have  $\|T(h_\alpha)\|=c\|h_\alpha\|$, that is,  $$\bigg[\mu(A)+|\alpha|^p\bigg(1-\frac{\epsilon}{p}\bigg)\mu(B)\bigg]^{\frac{1}{p}}=c[\mu(A)+|\alpha|^p\mu(B)]^{\frac{1}{p}},$$ which gives that $$c=\frac{\bigg[\mu(A)+|\alpha|^p\bigg(1-\frac{\epsilon}{p}\bigg)\mu(B)\bigg]^{\frac{1}{p}}}{[\mu(A)+|\alpha|^p\mu(B)]^{\frac{1}{p}}},$$
and this is a contradiction as $c$ depends on a scalar $\alpha$. Thus, the claim.

Now, we prove that $U_{\epsilon}$ preserves $\epsilon$-approximate B-J orthogonality. Let $0\neq f,\;g\in L^p(\mu,X)$ be such that $f\bj g$. In order to prove $U_{\epsilon}(f)\perp^{\epsilon}_{BJ}U_{\epsilon}(g)$, it is enough to prove that $|[U_{\epsilon}(g),U_{\epsilon}(f)]|\leq\epsilon\|U_{\epsilon}(f)\|\|U_{\epsilon}(g)\|$, by \cite[Proposition 3.1]{jacek}. Consider,
\begin{align*}
	\hspace{-1.3cm}
	[U_{\epsilon}(g),U_{\epsilon}(f)]&=\frac{1}{\|U_{\epsilon}(f)\|^{p-2}}\int\limits_{(Z(U_{\epsilon}(f)))^c}\|(U_\epsilon(f))(s)\|^{p-1}F_{(U_\epsilon(f))(s)}((U_\epsilon(g))(s))\,ds\\
	&=	\frac{1}{\|U_{\epsilon}(f)\|^{p-2}}\Bigg[\int\limits_{A\cap Z(f)^c}\|f(s)\|^{p-1}F_{f(s)}(g(s))\,ds + \\& \ \ \ \quad \quad \ \ \ \ \ \ \ \ \ \int\limits_{(S\setminus A)\cap Z(f)^c}\Bigg(1-\frac{\epsilon}{p}\Bigg)^{p-1}\|f(s)\|^{p-1}F_{f(s)}\bigg(\bigg(1-\frac{\epsilon}{p}\bigg)g(s)\bigg)\,ds\Bigg].
\end{align*}
Now, the space $L^p(\mu,X)$ is smooth, $X$ being smooth, so by \cite[Theorem 3]{giles}, $L^p(\mu,X)$ is a continuous semi-inner product space. This together with the fact that $f\bj g$ implies that  $[g,f]=\int\limits_{Z(f)^c}\|f(s)\|^{p-1}F_{f(s)}(g(s))\,ds=0$, \cite[Corollary 3.4]{jacek}. That is, $$\int\limits_{(S\setminus A)\cap Z(f)^c}\|f(s)\|^{p-1}F_{f(s)}(g(s))\,ds=-\int\limits_{A\cap Z(f)^c}\|f(s)\|^{p-1}F_{f(s)}(g(s))\,ds.$$ Thus, 
\begin{align*}
\hspace{-1.3cm}	|[U_{\epsilon}(g),U_{\epsilon}(f)]|&=\frac{1}{\|U_{\epsilon}(f)\|^{p-2}}\Bigg|\Bigg(1-\Bigg(1-\frac{\epsilon}{p}\Bigg)^{p}\Bigg)\int\limits_{A\cap Z(f)^c\cap Z(g)^c}\|f(s)\|^{p-1}F_{f(s)}(g(s))\,ds\Bigg|\\
	&\leq\frac{\epsilon}{\|U_{\epsilon}(f)\|^{p-2}}\int\limits_{A\cap Z(f)^c\cap Z(g)^c}\|f(s)\|^{p-1}\|g(s)\|\,ds\\
	&\leq\frac{\epsilon}{\|U_{\epsilon}(f)\|^{p-2}}\|g\vert_{(A\cap Z(f)^c\cap Z(g)^c)}\|\|f\vert_{(A\cap Z(f)^c\cap Z(g)^c)}\|^{p-1},
\end{align*}
 the first inequality follows from $\bigg[1-\bigg(1-\frac{\epsilon}{p}\bigg)^{p}\bigg]\leq\epsilon$ as $1<p<\infty$, and $\epsilon\in(0,1)$; and second follows from H\"{o}lder's inequality as done in the proof of \Cref{semiLp}. Also, observe that
\begin{align*}
	\|U_{\epsilon}(f)\|&=\Bigg[\int\limits_{A}\|f(s)\|^p\,ds+\int\limits_{S\setminus A}\Bigg(1-\frac{\epsilon}{p}\Bigg)^p\|f(s)\|^p\,ds\Bigg]^\frac{1}{p}\\
&	\geq\Bigg[\int\limits_{A}\|f(s)\|^p\,ds\Bigg]^\frac{1}{p}\\
&	\geq\Bigg[\int\limits_{A\cap Z(f)^c\cap Z(g)^c}\|f(s)\|^p\,ds\Bigg]^\frac{1}{p}\\
&	=\|f\arrowvert_{(A\cap Z(f)^c\cap Z(g)^c)}\|.
\end{align*}
Similarly, 	$\|U_{\epsilon}(g)\|\geq\|g\vert_{(A\cap Z(f)^c\cap Z(g)^c)}\|$, which further implies that $$|[U_{\epsilon}(g),U_{\epsilon}(f)]|\leq\frac{\epsilon}{\|U_{\epsilon}(f)\|^{p-2}}\|U_{\epsilon}(f)\|^{p-1}\|U_{\epsilon}(g)\|=\epsilon\|U_{\epsilon}(f)\|\|U_{\epsilon}(g)\|,$$ 
and this completes the proof.
\end{proof}

\end{document}